\documentclass{amsart}
\usepackage{amsmath,amssymb}
\usepackage[matrix,arrow,cmtip]{xy}
\renewcommand{\AA}{\mathbb{A}}
\newcommand{\CC}{\mathbb{C}}
\newcommand{\PP}{\mathbb{P}}
\newcommand{\ZZ}{\mathbb{Z}}
\newcommand{\QQ}{\mathbb{Q}}
\newcommand{\RR}{\mathbb{R}}
\newcommand{\FF}{\mathbb{F}}
\newcommand{\Aut}{\mathrm{Aut}}
\newcommand{\Gal}{\mathrm{Gal}}
\newcommand{\Div}{\mathrm{Div}}
\newcommand{\Pic}{\mathrm{Pic}}
\newcommand{\Jac}{\mathrm{Jac}}

\newcommand{\Exp}{\mathrm{Exp}}
\newcommand{\Log}{\mathrm{Log}}
\newcommand{\Spec}{\mathrm{Spec}}
\newcommand{\tng}{\mathrm{Tang}}

\newcommand{\Prym}{\mathrm{Prym}}

\newcommand{\id}{\mathrm{id}}
\newcommand{\Res}{\mathrm{Res}}

\newcommand{\rk}{\mathrm{rk}}

\newcommand{\Kbar}{{\overline{K}}}

\newcommand{\act}[2]{{}^{#1}\!#2}
\newcommand{\fp}{\mathfrak{p}}

\newcommand{\fD}{\mathfrak{D}}

\newcommand{\sD}{\mathcal{D}}

\newcommand{\sK}{\mathcal{K}}

\newcommand{\sO}{\mathcal{O}}

\newcommand{\vx}{\mathbf{x}}
\newlength{\underscorelength}
\newcommand{\deadus}{\settowidth{\underscorelength}{{\rm \_}}\hspace{-0.5\underscorelength}\_\hspace{-0.5\underscorelength}}
\newcommand{\Sha}{{\rm I{\deadus}I{\deadus}I}}

\renewcommand{\phi}{\varphi}
\newtheorem{theorem}{Theorem}[section]
\newtheorem{lemma}[theorem]{Lemma}
\newtheorem{corollary}[theorem]{Corollary}
\newtheorem{prop}[theorem]{Proposition}
\theoremstyle{definition}
\newtheorem{definition}[theorem]{Definition}

\newtheorem{remark}[theorem]{Remark}
\begin{document}
\bibliographystyle{plain}
\title{The arithmetic of Prym varieties in genus $3$}
\parskip=0.1pt
\author{Nils Bruin}
\address{Department of Mathematics, Simon Fraser University,
Burnaby BC, Canada V5A 1S6}
\email{bruin@member.ams.org}
\thanks{The research in this paper is partially supported by NSERC}


\date{May 30, 2006}

\begin{abstract}
Given a curve of genus $3$ with an unramified double cover, we give an explicit
description of the associated Prym variety. We also describe how an unramified
double cover of a non-hyperelliptic genus $3$ curve can be mapped into
the Jacobian of a curve
of genus 2 over its field of definition and how this can be used to do
Chabauty- and Brauer-Manin type calculations for curves of genus 5 with an
unramified involution. As an application, we determine the rational points on a
smooth plane quartic with no particular geometric properties and give
examples of curves of genus $3$ and $5$ violating the Hasse-principle. We also show
how these constructions can be used to design smooth plane quartics with
specific arithmetic properties. As an example, we give a smooth plane quartic
with all $28$ bitangents defined over $\QQ(t)$, By specialization, this also
gives examples over $\QQ$.
\end{abstract}
\maketitle
\section{Introduction}
In this article, we investigate the arithmetic of unramified double covers of
non-hyperelliptic curves of genus $3$. This research is inspired by the recent
success in applying the theory of unramified double covers to hyperellipitic
curves to the problem of determining the set of rational points on such curves
(\cite{weth:phdthesis}, \cite{bruin:tract}, \cite{bruin:crelle},
\cite{brufly:tow2cov}).
Combined with explicit Chabauty-methods, this has yielded very practical methods
for obtaining often sharp bounds on the number of rational points on a
hyperelliptic curve. In the hyperelliptic case, the construction of unramified
covers is particularly easy to make explicit using Kummer theory.
There has been some limited work investigating how these ideas may
be generalised to unramified covers of higher degree of hyperelliptic curves
(\cite{brufly:ncov}).

In this article, we generalise in a different direction. We derive the general
form of a curve $C$ of genus $3$ over a field $K$ of characteristic $0$ that allows an
unramified double cover $\pi: D\to C$ defined over $K$. A non-hyperelliptic
curve of genus $3$ with an unramified double cover allows a smooth plane
projective model of the form
$$C: Q_1(u,v,w)Q_3(u,v,w)=Q_2(u,v,w)^2$$
where $Q_1,Q_2,Q_3\in K[u,v,w]$ are quadratic forms.
The unramified double cover has a canonical model in $\PP^5$ of the form
$$D:\left\{
\begin{array}{rcl}
Q_1(u,v,w)&=&r^2\\
Q_2(u,v,w)&=&rs\\
Q_3(u,v,w)&=&s^2\;.
\end{array}\right.$$

We also describe how the
associated Prym variety (the connected component containing $0$ of the kernel of
$\pi_*:\Jac(D)\to\Jac(C)$) can be described explicitly as $\Jac(F)$ of some
genus $2$ curve $F$ over $K$:
$$F:y^2=-\det(Q_1+2xQ_2+x^2 Q_3)$$
where we regard $Q_i$ as the symmetric $3\times 3$ matrices corresponding to
the quadratic forms they represent. Combined with earlier work, this gives a
complete description of how principally polarised Abelian surfaces arise as
Prym varieties in genus $3$ over non-algebraically closed base fields of
characteristic $0$ (see Theorem~\ref{thrm:prymclass}).

The description of $\Prym(D/C)$ as $\Jac(F)$ over $\CC$
can already be found in \cite[Exercise~VI.F]{ACGH:curves}. Alternatively, if $C$
has a rational point, then the \emph{trigonal construction} gives a
Galois-theoretic way of describing the Prym variety as the Jacobian of some
curve (see \cite{donagi:fibprym} and \cite{recillas:trigonal}). Here, we refine the
description to arbitrary characteristic $0$ base fields. The results should
generalise to arbitrary odd characteristic.

We show how $D$ can be mapped into $\Prym(D/C)$ over $K$, without additional
assumptions on $D$ (see Proposition~\ref{prop:Demb}). We can thus use
Chabauty-like methods on genus $5$ curves
in Abelian surfaces to get an in practice often sharp bound on the number of rational points on $D$.
This can in turn be used to get a corresponding bound on the number of rational
points on $C$. As an example, we prove
\begin{prop}\label{prop:ex1}
For the curve
$$C : (4u^2- 4vw + 4w^2)(2u^2 + 4uv + 4v^2) =( 2u^2 + 2uw- 4vw + 2w^2)^2$$
we have $C(\QQ)=\{(0:1:0)\}$.
\end{prop}
The proof avoids computing the Mordell-Weil group of $\Jac(C)$. Instead, it uses
the computationally more accessible Mordell-Weil group of an Abelian surface.

Furthermore, the method does not put restrictions on the geometry of $C$: Over
some field extension $L$, the group scheme $\Jac(C)[2]$ has a non-trivial
rational point and therefore an unramified cover $D$ of the desired type and a
map $D\to\Prym(D/C_L)$.
In complete analogy to the treatment of genus $2$ in \cite{bruin:crelle}, by
taking the Weil-restriction of scalars, we obtain
$\Re_{L/K}(D)\to\Re_{L/K}(\Prym(D/C_L))$. Inside of $\Re_{L/K}(D)$, we can find a
curve $\tilde{D}$ corresponding to $\Re_{C_L/C_K}(D/C_L)$. We can then apply
Chabauty-like methods to $\tilde{D}\to\Re_{L/K}(\Prym(D/C_L))$ to determine
$\tilde{D}(K)$ and $C(K)$.
The group $\Re_{L/K}(\Prym(D/C_L))(K)\simeq \Prym(D/C_L)(L)$ would be the
hardest ingredient to obtain and it should be noted that the computations
involved would probably be prohibitive, except for very low degree $L$.

Since the mapping of $D$ into $\Prym(D/C)$ does not require a rational point on
$D$, we can apply this construction to prove that $D(\QQ)$ is empty, even if $D$
does have points everywhere locally. It is reassuring that, at least
conjecturally, our computations  correspond to determining part of the
Brauer-Manin obstruction of $D$:

Suppose that $D$ is a curve of genus larger than $1$, defined over a number
field $K$, with points everywhere locally. Assuming $\Sha(\Jac(D)/K)$ is finite, a failure
for $D$ to have a $K$-rational degree $1$ divisor class would be due to the
Brauer-Manin obstruction (see \cite[Corollary 6.2.5]{skor:tor}). Otherwise,
using an Abel-Jacobi embedding, $D$ can be considered a subvariety of $\Jac(D)$.
Scharaschkin \cite{schar:phdthesis} proves that if $D(\AA_K)$ misses the
topological closure of $\Jac(D)(K)$ in $\Jac(D)(\AA_K)$, then this is due to the
Brauer-Manin obstruction on $D$ if $\Sha(\Jac(D)/K)$ is finite. 

In our computations, we show that the image of $D(\AA_{\QQ})$ in
$\Prym(D/C)(\AA_{\QQ})$ misses the closure of $\Prym(D/C)(\QQ)$ by combining the
local information at a finite number of primes.

As an example, we prove
\begin{prop}\label{prop:ex2}
The genus $5$ curve
$$D:\left\{
\begin{array}{rcl}
(v^2 + vw - w^2)&=&r^2\\
(u^2 - v^2 - w^2)&=&rs\\
(uv + w^2)&=&s^2
\end{array}\right.$$
and the genus $3$ curve
$$C: (v^2 + vw - w^2)(uv + w^2)=(u^2 - v^2 - w^2)^2$$
both have points everywhere locally, but they have no rational points.
\end{prop}
For this example, we find $\Prym(D/C)(\QQ)\simeq \ZZ\times\ZZ$. This illustrates
that the method used is not a special case of a Chabauty-type argument.

Additionally, we investigate the arithmetic implications of the geometric
description of the fibres of the Prym map between moduli spaces, corresponding
to the functor $\{\pi:D\to C\}\mapsto\{\Prym(D/C)\}$,
given in \cite{verra:prym}. For a general Abelian surface $A$,
the Kummer surface
$\sK=A/\langle \pm 1\rangle$ has a singular quartic model in $\PP^3$, with a
singular locus consisting of $16$ points, corresponding to $A[2]$. Translation by $A[2]$
induces automorphisms on $\sK$, given by linear transformations on $\PP^3$. We
write $\widehat{\PP}^3$ for the space of planes in $\PP^3$.
Verra \cite{verra:prym} shows that over $\CC$, the fibre of the Prym map over a principally
polarised surface $A$, is birational
to $\widehat{\PP}^3/A[2]$. In fact, he gives a very precise description of the
fibre as a blow-up of this space, where the exceptional components contain the
moduli points corresponding to hyperelliptic or otherwise degenerate curves $C$.

In particular, a non-hyperelliptic genus $3$ curve $C$ over $\CC$ which has a
double cover $D$ such that $\Prym(D/C)=A$ can be obtained as a plane section of
$\sK$, with $D$ the pull-back of $C$ to $A$. Any two such plane sections of
$\sK$ in the same $A[2]$-orbit give isomorphic covers
$D/C$.

In this article we explain how, given a genus $2$ curve $F$ over a number field
$K$ and a genus $3$ plane section $C$ of $\sK=\Jac(F)/\langle \pm 1\rangle$, we
can obtain a model for $C$ of the type $Q_1Q_3=(Q_2)^2$. Since a sufficiently general plane
section of $\sK$ is non-singular and thus of genus $3$, it shows that any
Jacobian of a genus $2$ curve over $K$ can be realised as a Prym variety of a
non-hyperelliptic genus $3$ curve \emph{over $K$}, or stated more amusingly in
an elementary fashion:
\begin{prop}\label{prop:fdet}
Let $f\in\QQ[x]$ be a square-free polynomial of degree $5$ or $6$. Then there
exist symmetric matrices $M_1,M_2,M_3\in\QQ^{3\times 3}$ such that
$$f = \det(M_1+x M_2+ x^2 M_3).$$
\end{prop}
We use this construction to obtain a systematic way of constructing curves of
genus $3$ with all $28$ bitangents defined over a non-algebraically closed
field, for instance $\QQ(t)$ (see
Section~\ref{sec:fibprym}). By specialisation of $t$, this gives infinitely many
examples over $\QQ$. This strenghtens a
result in \cite{edge:realbitang}, where an example is given with all bitangents
defined over $\RR$. See \cite{scha:selfun} or \cite{cayley:bitang} for
an approach via interpolation.

Finally, it should be noted that not all covers $\pi:D\to C$ defined over $K$
with $\Prym(D/C)=A$ and $C$ non-hyperelliptic have $C$ occurring as a plane
section of the associated Kummer surface $\sK$ defined over $K$, since not all
rational points of $\widehat{\PP}^3/A[2]$ are covered by rational points of
$\widehat{\PP}^3$. In particular, the $16$ tropes of $\sK$ give rise to a
Galois-stable set of $16$ bitangents to $C$. Not all smooth plane quartics
have such a configuration of bitangents.

\section{Preliminaries}

First we fix some notation. Let $K$ be a field and let $C$ be a complete,
absolutely irreducible algebraic curve over $K$. We write $\kappa_C$ for a
canonical divisor on $C$. For a divisor $\fD\in\Div(C)$, we write $[\fD]$ for
its class in $\Pic(C)$ and $|\fD|$ for the complete linear system corresponding
to $\fD$ and $l(\fD)=\dim |\fD|$. We say that a divisor $\fD$ on $C$
is a $g^r_d$ if $\deg(\fD)=d$ and $l(\fD)> r$. We write
$W^r_d\subset\Pic^d(C)$ for the classes of $g^r_d$s. By abuse of notation, we
will also write $W^r_d$ for the corresponding subscheme of the scheme
representing the functor $\Pic^d(C)$.

\subsection{Unramified double covers and Prym varieties}

Let $\pi: D\to C$ be an unramified finite morphism of degree $2$ between curves
over a field $K$ and let $\iota:D\to D$ be the involution of $D$ over $C$.
It follows by Riemann-Hurwitz that $g(C)>0$ and that
$$g(D)=2g(C)-1.$$
The kernel of $\pi_*:\Jac(D)\to\Jac(C)$ has two connected components. The
component that contains $0$ coincides with the image of
$(\id_*-\iota_*):\Jac(D)\to\Jac(D)$.

\begin{definition}
Let $\pi: D\to C$ be an unramified finite morphism of degree $2$ between curves
over a field $K$. We write $\Prym(\pi)=\Prym(D/C)$ for the connected component
of the identity-element of $\ker(\pi_*:\Jac(D)\to\Jac(C))$. We call this the
\emph{Prym variety} of $D/C$.
\end{definition}

Thus, $\Prym(D/C)$ is an Abelian subvariety of $\Jac(D)$. The principal
polarisation on $\Jac(D)$ restricts to a principal polarisation on $\Prym(D/C)$.
Historically, Prym varieties were considered interesting primarily because they
give examples of principally polarised Abelian varieties that are not Jacobian
varieties. However, if $\dim(\Prym(D/C))\leq 2$, then $\Prym(D/C)$ generally is
a Jacobian variety\footnote{It can also be the product of two elliptic curves,
in which case it is a generalised Jacobian.}.  See \cite[VI-C]{ACGH:curves} or
\cite{mum:prym1} for details.

\subsection{Prym varieties in the hyperelliptic case}
\label{sec:prymhyp}

Contrary to the general situation,
the Prym variety associated to an unramified double cover of a hyperelliptic
curve is closely related to a Jacobian variety. In fact, the Prym variety is
isomorphic to the product of Jacobian varieties of subcovers, which themselves
are again hyperelliptic.

Let us first assume that $C$ is a double cover of a $\PP^1$. Then $C$ has an
affine model of the form
$$C: y^2=f(x)$$
where $f\in K[x]$ is a square-free polynomial of degree $2g(C)+2$. Kummer theory
tells us exactly what the unramified degree $2$ extensions of $K(C)$ are. For
any factorisation $f=f_1f_2$, with $f_1,f_2\in K[x]$ and of even degree, we have
a curve $D$ given by the affine model
$$\left\{\begin{array}{rcl}
y_1^2&=&f_1(x)\\
y_2^2&=&f_2(x)\\
\end{array}\right.$$
and an unramified morphism of degree $2$
$$\begin{array}{cccccc}
\pi&:&D&\to&C&\\
&&(x,y_1,y_2)&\mapsto&(x,y_1y_2)&=(x,y).
\end{array}$$
Then there are the two obvious curves
$$\begin{array}{cccc}
F_1:&y_1^2&=&f_1(x)\\
F_2:&y_2^2&=&f_2(x).
\end{array}$$
with the obvious projections $\pi_1:D\to F_1$ and $\pi_2:D\to F_2$. This yields
the familiar diagram associated to biquadratic extensions.
$$\xymatrix{
&D\ar[dr]^{\pi_2}\ar[d]_\pi\ar[dl]_{\pi_1}\\
F_1\ar[dr]_{x}&C\ar[d]_x&F_2\ar[dl]^{x}\\
&\PP^1
}$$

\begin{prop}
\label{prop:hypprymsplit}
Let $C,D,F_1,F_2,\pi,\pi_1,\pi_2$ be defined as above. Then
$$\pi_1^*\times\pi_2^*:\Jac(F_1)\times\Jac(F_2)\to\Prym(D/C)$$
is an isomorphism of Abelian varieties.
\end{prop}

\begin{proof}
First, we prove that $\pi_1^*$ indeed maps $\Jac(F_1)$ into $\Prym(D/C)$.
To that end, take the generic point $(x_1,y_1)\in F_1$. We have
$\pi_1^*(x_1,y_1)=(x_1,y_1,\sqrt{f_2(x_1)})+(x_1,y_1,-\sqrt{f_2(x_1)})$.
Under $\pi$, this maps to the divisor cut out by $x=x_1$. This shows that
$\pi_*\pi_1^*:\Jac(F_1)\to\Jac(C)$ is constant and hence
$(\pi_*\pi_1^*)|_{\Jac(F_1)}$ is the zero map. By symmetry it follows that
$\Jac(F_1)\times\Jac(F_2)$ lands in the kernel of $\pi_*$ under
$\pi_1^*+\pi_2^*$ and, since it is connected and covers $0\in\Jac(D)$, lands in
$\Prym(D/C)$.

The fact that $(\pi_i)_*\circ\pi_i^*=2|_{\Jac(F_i)}$ for $i=1,2$ already
assures that $\Jac(F_1)\times\Jac(F_2)$ is isogenous to $\Prym(D/C)$. A quick
way to see that they are actually isomorphic is by noting that $F_1$ and $F_2$
can be arbitrary hyperelliptic curves and that, by construction, the isogeny
would have to depend functorially on $F_1$ and $F_2$. In general,
$\Jac(F_1)\times\Jac(F_2)$ has no non-trivial polarisation-preserving
isogenies and hence there are no other candidates for $\Prym(D/C)$.
\end{proof}

\begin{remark}
There are more complicated situations than the split situation described in
\ref{prop:hypprymsplit}. For instance, if $L$ is a quadratic extension of $K$
with conjugation $\sigma:L\to L$ over $K$ and $f_2=\act{\sigma}{f_1}$, then
$F_2=\act{\sigma}{F_1}$. There is still an unramified double cover $D/C$ over
$K$ associated with this splitting. The corresponding Prym variety is the
Weil restriction $\Re_{L/K}\Jac(F_1)$.
\end{remark}

\begin{remark}
In general, a hyperelliptic curve $C$ over $K$ is a double cover of a curve $L$
of genus $0$. We can express $K(L)$ as some quadratic extension of $K(\PP^1)$.
Relative Kummer theory allows us to describe $\Prym(D/C)$ in terms of Jacobians
of subcovers of $D/L$ in exactly the same way as above.
\end{remark}

\subsection{Linear subspaces on quadrics}
\label{sec:ruling}

It is well known that on a non-singular quadric in $\PP^3$, there are two
rulings of lines with the property that a line from one ruling intersects a
unique line from the opposite ruling and that two lines from the same ruling do
not intersect. We will need a simple lemma that classifies whether the two
rulings are split or quadratic conjugates.

\begin{lemma}
\label{lemma:quadsplit}
Let $K$ be a field of characteristic different from $2$ and let
$M\in K^{4\times 4}$ be a symmetric matrix describing a non-singular quadric
$Q\subset\PP^3$. The two rulings of lines on $Q$ are individually defined over
$K$ exactly if $\det(M)$ is a square in $K$. Otherwise, they are quadratic
conjugate.
\end{lemma}

\begin{proof}
First assume we have a point $\vx_0\in Q(K)$. Let
$V=\{\vx:\act{t}{\vx_0}M\vx=0\}$ be the plane tangent to $Q$ at $\vx$. The
intersection $V\cap Q$ consists exactly of two lines through $\vx_0$, one from
each of the rulings. With a change of basis, we can assume that
$\vx_0=(1:0:0:0)$ and that the points $\vx_1=(0:1:0:0)$ and $\vx_2=(0:0:1:0)$
lie on $V$ as well. It follows that
$$M=\left(
\begin{array}{cccc}
0&0&0&d\\
0&a&b&*\\
0&b&c&*\\
d&*&*&*\\
\end{array}
\right);\quad
\det(M)=d^2(b^2-ac)$$
We see that $Q$ intersected with the line through $\vx_1,\vx_2$ is described by
the equation $a x_1^2+2b x_1x_2+c x_2^2=0$, which is split exactly if $b^2-ac$
is a square. The lemma follows.

If $Q(K)$ is empty, then we base change to $K(Q)$, where we have the generic
point $(x_0:x_1:x_2:x_3)\in Q(K(Q))$. Since $K$ is algebraically closed in
$K(Q)$, the pair of rulings (defined over $K$) is split over $K(Q)$ if and only
if they are split over $K$. Furthermore, $\det(M)$ is a square in $K$ if and
only if it is in $K(Q)$.
\end{proof}

\section{Non-hyperelliptic curves of genus $5$ with an unramified involution}
\label{sec:genus5}

Let $K$ be a field of characteristic $0$ and let $D$ be a non-hyperelliptic
curve of genus $5$ with an unramified involution $\iota$. Let
$\pi:D\to C=D/\langle \iota\rangle$ be the quotient map by the action
of $\iota$. The Riemann-Hurwitz formula yields that $C$ is of genus $3$.
Let $\kappa_C$ be a canonical divisor on $C$ and let $\langle
u,v,w\rangle=|\kappa_C|$ be coordinates on the associated canonical model of
$C$. Note that we do not insist that $C$ is non-hyperelliptic. By abuse of
notation, we also write $u,v,w$ for the pull-backs along $\pi$ of the
corresponding functions on $C$.
We write $\kappa_D=\pi^*\kappa_C$. There are functions $r,s$ on $D$ with
$r+r\circ\iota=s+s\circ\iota=0$ such that
$\langle u,v,w,r,s\rangle=|\kappa_D|$.

We identify $D$ with the canonical model associated to $\kappa_D$, being the
image of $(u:v:w:r:s): D\to\PP^4$. In this notation,
$$\iota:(u:v:w:r:s)\mapsto(u:v:w:-r:-s).$$

We let $\Lambda$ be the linear system of quadrics containing $D\subset\PP^4$. A
simple comparison of the dimensions of $l(\kappa_D)=5$ and $l(2\kappa_D)=12$
yields that $\Lambda\simeq\PP^2$. Let $Q_1,Q_2,Q_3\in\Lambda(K)$ be quadrics
generating $\Lambda$ and let $(\lambda_1:\lambda_2:\lambda_3)$ be the
corresponding coordinates on $\Lambda$.

Note that $|\kappa_D|$ has a decomposition in the $+1$-eigenspace $\langle
u,v,w\rangle$ and the $-1$-eigenspace $\langle r,s\rangle$ of $\iota$. The
involution $\iota$ acts identically on the
corresponding linear subspaces $\{r=s=0\},\{u=v=w=0\}\subset\PP^4$.

\begin{lemma}\label{lemma:nontrig}
With the notation above, $D$ is not a trigonal curve.
\end{lemma}

\begin{proof}
Since a trigonal curve remains trigonal upon base extension, it suffices to
prove the lemma for algebraically closed $K$.

We argue following \cite[p.~207]{ACGH:curves}.
Suppose that $\fD\in\Div(D)$ is a $g^1_3$. Then $\kappa_D-\fD$ is a $g^2_5$, so
$D$ has a plane quintic model. This model must have a unique singularity $P_0$
and the divisors in $|\fD|$ are cut out by lines through $P_0$. 

As is argued in \cite[p.~207]{ACGH:curves} (using the Brill-Noether Residue
Theorem, see for instance \cite[Ch.~8]{fulton:curves}), a curve of genus $5$
can
have at most $1$ divisor class of type $g^1_3$. Assume that $D$ is trigonal and
let $\sD$ be the divisor class of type $g^1_3$. Then $\iota_*:\Div(D)\to\Div(D)$
induces an involution on $|\sD|\simeq\PP^1$. Let $\phi=\phi_\sD: D\to\PP^1$ be
the trigonal map on $D$. Then $D\to |\sD|$ defined by $P\mapsto
\phi^*(\phi(P))$ shows that $|\sD|$ is naturally the $\PP^1$ of which $D$ is
the degree $3$ cover.

Since $|\sD|$ contains all effective $g^1_3$s on $D$, the restriction of
$\iota^*$ to $|\sD|$ yields an involution on $|\sD|$. Since involutions on
$\PP^1$ have two fixed points, there are two
effective $g^1_3$s on $D$ on which are fixed under $\iota_*$. It follows that
$\iota$ permutes the support of one such divisor, which consists of $3$ not
necessarily distinct points. Since $\iota^2=1$, it follows that at least one
point is fixed under $\iota$ and hence that $\pi$ has ramification.
\end{proof}

\begin{lemma}\label{lemma:quadsect}
With the notation above, $D=Q_1\cap Q_2\cap Q_3$. Furthermore, the $Q_i$ can be
chosen to be non-singular and $D$ misses the singular locus of any quadric
containing $D$.
\end{lemma}

\begin{proof}
First note that if the statement of the lemma is false, then it is also false
over the algebraic closure of $K$. Therefore, it suffices to prove the lemma for
algebraically closed $K$.
By Petri's Theorem \cite[p.~131]{ACGH:curves}, a canonical model of a
non-hyperelliptic, non-trigonal curve of genus $5$ is the intersection of
quadrics. 

Next we show that a quadric $Q\in\Lambda$ cannot be singular at $D$.
Suppose that $P_0\in D$ is a singular point of some quadric $Q\in\Lambda$. Note
that, if $\rk Q<3$, then there is a $\PP^3\subset Q$. The space $\Lambda$
restricted to that $\PP^3$ would be at most a pencil of quadrics and hence
contain an intersection of $2$ quadrics. This would imply that $D$ has a
component of genus at most $1$, which contradicts that $D$ is a curve of
genus $5$.

Hence, $Q$ is of rank $3$ or $4$, which implies that $Q$ contains a $\PP^1$ of
planes through the singular point $P_0\in D$. On each such plane $V$, the
restriction of $\Lambda$ is a pencil of conics and hence has a base locus of
degree $4$. Since $P_0\in V\cap D$ is contained in
that base locus, this realises $D$ as a degree $3$ cover of that $\PP^1$. This
contradicts that $D$ is non-trigonal.

By Bertini's Theorem \cite[III,~10.9.2]{har:AG} it follows that a
general member of $\Lambda$ is non-singular and hence we can choose the $Q_i$
to be non-singular.
\end{proof}

\begin{lemma}\label{lemma:decomp}
With the notation above, we have that $\iota$ acts trivially on
$\Lambda$. Equivalently, for $Q\in\Lambda$, we can find quadratic forms
$Q^+\in K[u,v,w]$ and $Q^-\in K[r,s]$ such that $Q$ is given by the equation
$$Q^+(u,v,w)+Q^-(r,s)=0.$$
\end{lemma}

\begin{proof}
First note that $\iota$ preserves $D$ and hence preserves $\Lambda$. Suppose
there is $Q\in\Lambda(\Kbar)$ such that $Q\neq\iota(Q)$. Then $Q-\iota(Q)\in\Lambda(\Kbar)$ is
not $0$. On the other hand, since $\iota$ acts trivially on $\{r=s=0\}$, we have
that $Q-\iota(Q)$ restricted to $\{r=s=0\}$ is zero. Hence, $\Lambda$ restricted
to $\{r=s=0\}$ is a pencil of plane conics and has a base locus of degree $4$.
It follows that $\iota$ is ramified, which contradicts the assumptions.

For an equation of any quadric in $\Lambda(\Kbar)$, this means that the monomials
$ur,vr,\ldots,ws$ cannot occur.
\end{proof}

\section{Special divisor classes of degree $4$}
\label{sec:W14}

As in the previous section, let $D$ be a non-hyperelliptic curve of genus $5$
over a field $K$ of characteristic $0$ with an unramified involution
$\iota:D\to D$. We adopt the other notation from the previous section as well.
We consider the scheme of special divisors
$$W^1_4(D)=\{\sD\in\Pic^4(D):l(\sD)\geq 2\}.$$
From the Riemann-Roch formula it follows that the residuation map
$\sD\mapsto[\kappa_D]-\sD$ defines an involution on $W^1_4$. 

We denote the locus of singular quadrics in $\Lambda$ by
$$\Gamma=\{Q\in\Lambda:\det(Q)=0\}.$$
By Lemma~\ref{lemma:decomp} we have a decomposition
$\Gamma=\Gamma^+\cup\Gamma^-$, with equations
$$\begin{array}{rrcl}
\Gamma^+:& \det(\lambda_1 Q_1^+ +\lambda_2 Q_2^+ + \lambda_3 Q_3^+)&=&0,\\
\Gamma^-:& \det(\lambda_1 Q_1^- +\lambda_2 Q_2^- + \lambda_3 Q_3^-)&=&0.
\end{array}$$
By Lemma~\ref{lemma:quadsect}, $\Gamma$ is $1$-dimensional. In fact, it is
straightforward to check that
$$\Gamma'=\{Q\in\Lambda: \rk Q = 3\}$$
is $0$-dimensional and is the singular locus of $\Gamma$.

\begin{lemma}
\label{lemma:quadplanes}
Let $\fD,\fD'\in\Div(D)$ be effective divisors of degree $4$ with $l(\fD)=2$.
Then the following hold.
\begin{itemize}
\item[(i)] There is a unique $2$-plane $V_\fD$ such that $\fD=D\cdot V_\fD$.
\item[(ii)] There is a unique quadric $Q_\fD\in\Lambda$ which vanishes on
$V_\fD$. In fact, $Q_\fD\in\Gamma$.
\item[(iii)] If $V_\fD$ and $V_{\fD'}$ meet in a line, then
$[\fD+\fD']=[\kappa_D]$.
\item[(iv)] If $Q_\fD = Q_{\fD'}$ then $[\fD'+\fD]=[\kappa_D]$ or $[\fD-\fD']=0$.
\item[(v)] If $Q_\fD \neq Q_{\fD'}$ then $[\fD'+\fD]\neq[\kappa_D]$ and
$[\fD-\fD']\neq 0$.
\end{itemize}
\end{lemma}
\begin{proof}
\noindent (i): The geometric formulation of the Riemann-Roch Theorem
\cite[p.\ 12]{ACGH:curves} states that, for a divisor
$P_1+\cdots+P_r$ with $r\leq g$, we have $$l(P_1+\cdots+P_r)=r+1-\rk\langle
P_1,\ldots,P_r\rangle,$$ so one can take $V_\fD$ to be the plane spanned by the
support of $\fD$ over $\Kbar$.

\medskip
\noindent (ii): Since the restriction of $\Lambda$ to $V_\fD$ has a base locus of
degree $4$, it is a pencil of conics. Hence there is a unique quadric
$Q_\fD\in\Lambda$ that vanishes on $V_\fD$. Since a quadric in $\PP^4$
containing a $2$-plane is necessarily singular, it follows that
$Q_\fD\in\Gamma$.

\medskip
\noindent (iii): Two $2$-planes meeting in a line lie in a $3$-plane, say $W$.
Since $\fD+\fD'$ is effective of degree $8$ it equals the
hyperplane section $D\cdot W$. Hyperplane sections of canonical models are
canonical divisors.

\medskip

\noindent (iv): First suppose that $Q_\fD$ is of rank $4$. Then $Q_\fD$ is a
cone over a nonsingular quadric in $\PP^3$. The two line rulings on that quadric
give rise to two plane ``rulings'' on $Q_\fD$. Planes in opposite rulings meet
in a line and planes in the same ruling meet only in the singular point of
$Q_\fD$. Hence, if $V_\fD$ and $V_{\fD'}$ belong to opposite rulings, then by
(iii) we have $[\fD+\fD']=[\kappa_D]$. If $V_\fD$ and $V_{\fD'}$ belong to the
same ruling then both $\fD$ and $\fD'$ are residual to an arbitrary divisor from the
opposite ruling and hence linearly equivalent.

If $Q_\fD$ is of rank $3$ then $Q_\fD$ is a cone over a singular quadric in
$\PP^3$, so there is one ruling of planes on $Q_\fD$ and any two of these
meet in a line. It follows by (iii) that $[\fD+\fD']=[\kappa_D]$. Since
this holds for any $\fD'$ with $Q_{\fD'}=Q_\fD$, it follows that that
$[\fD-\fD']=0$.

\medskip

\noindent (v): Given $V_\fD$, there is a map
$$\PP^1\to\{\mbox{planes in }Q_\fD\}$$
parametrising the ruling on $Q_\fD$ containing $V_\fD$. Each of the planes $V$
in the image of this map cuts out an effective divisor $\fD'$ equivalent to
$\fD$. By construction, $Q_{\fD'}=Q_\fD$. This gives an explicit realisation of
$\PP^1\simeq |\fD|$, so this accounts for all divisors linearly equivalent to
$\fD$.

\end{proof}

\begin{corollary}
The map
$$\begin{array}{ccc}
W^1_4(D)&\to&\Gamma\\
{}[\fD]&\mapsto& Q_\fD
\end{array}$$
is well-defined and realises $W^1_4(D)$ as a double cover of $\Gamma$, ramified
over $\Gamma'$. The involution of $W^1_4(D)$ over $\Gamma$ corresponds to the
residuation map $\sD\mapsto[\kappa_D]-\sD$.
\end{corollary}

\section{The Prym variety in genus $3$.}
\label{sec:prym}
We write $F$ for the union of components of $W^1_4(D)$ above $\Gamma^-$. If $\sD\in F$
then $Q_\sD$ has a singular point on $\{u=v=w=0\}$. Since the plane $V_\sD$
passes through this singularity, we see that $\pi(V_\sD)$ is a line and hence
that $\pi_*(\sD)=[\kappa_C]$. It follows that the map
$$\begin{array}{ccc}
  F\times F&\to& \Pic^0(D)\\
  \sD_1,\sD_2&\mapsto&\sD_1+\sD_2-[\kappa_D]
\end{array}$$
maps $F\times F$ into $\ker(\pi_*)$.

Depending on the type of $\Gamma^-$, this gives us different descriptions of
$\Prym(D/C)$. The case numbering is in correspondence with
Table~\ref{tab:prymstruct}.

\medskip
\noindent\emph{Case 0}: $\Gamma^-$ is a double counting line. This case does not
occur because assuming it does leads to a contradiction.
By choosing coordinates appropriately,
$\Gamma^-$ is described by the equation $\lambda_1^2=0$. It follows that
$Q_1^-+\lambda_2Q_2^-+\lambda_3 Q_3^-$ is non-singular for all
$\lambda_2,\lambda_3$. It follows that $Q_2^-=Q_3^-=0$. Therefore, $D$ has an
intersection with $\{u=v=w=0\}$ and thus $\iota$ would be ramified.

\medskip
\noindent\emph{Case 2}: $\Gamma^-$ is a split singular conic, i.e.,
$L_1\cup L_2$. In that situation, $F$ consists of two
components $E_1$ and $E_2$, covering $L_1$ and $L_2$ respectively. Each of these
 covers is ramified at a degree $4$ locus: the intersection of $L_i$ with the
 other components of $\Gamma$. Hence $E_1$ and $E_2$ are curves of genus $1$. In
 fact, each has a rational point above $L_1\cap L_2$, so we can identify them
 with their Jacobians. We
 have the map
$$\begin{array}{rcl}
E_1 \times E_2&\to&\Jac(D)\\
(\sD_1,\sD_2)&\mapsto&\sD_1+\sD_2-[\kappa_D]
\end{array}.$$
Note that $(\sD_1,\sD_2)$ only map to $0$ if $Q_{\sD_1}=Q_{\sD_2}$. This can
only happen above the (ramified) point $L_1\cap L_2$, so the map is an
injection. Since $E_1\times E_2$ is connected, the image is contained in
$\Prym(D/C)$ and because $E_1\times E_2$ is an Abelian surface itself, we have equality:
$$\Prym(D/C)\cong E_1\times E_2.$$

\medskip
\noindent\emph{Case 3}: $\Gamma^-$ is a non-split singular conic. Then, over some
quadratic extension $K(\sqrt{d})$ of $K$, the conic $\Gamma^-$ splits
and the analysis above applies. It follows that in that situation $E_1$ and
$E_2$ are elliptic curves that are conjugate with respect to $K(\sqrt{d})/K$. By
Weil-restriction, it follows that
$$\Prym(D/C)\cong \Re_{K(\sqrt{d})/K}(E_1).$$

\medskip
\noindent\emph{Case 4}: $\Gamma^-$ is a non-singular conic. In that case,
$Q_1^-,Q_2^-,Q_3^-$ are $K$-linearly independent and therefore span the space of
quadratic forms in $r,s$. Without loss of generality, we can assume
$$\begin{array}{rrcl}
Q_1:&Q_1^+(u,v,w)=r^2,\\
Q_2:&Q_2^+(u,v,w)=rs,\\
Q_3:&Q_3^+(u,v,w)=s^2.
\end{array}$$
It follows that $\Gamma^-$ is given by the equation
$4\lambda_1\lambda_3=\lambda_2^2$ and we have a parametrisation
$$\begin{array}{ccc}
 \PP^1&\to&\Gamma^-\\
 (x:1)&\mapsto&(1:2x:x^2).
\end{array}$$
The curve $F$ is a double cover of the $\Gamma^-$, ramified above
$\Gamma^+\cap\Gamma^-$. Using the parametrisation above we get an equation of
the form
$$\det(Q_1^++2xQ_2^++x^2Q_3^+)=\delta y^2.$$
for some $\delta\in K^*$. In order to determine the correct value of $\delta$,
suppose we have $\sD\in F$ such that $2\sD\neq [\kappa_D]$. Then $|\sD|$ is
cut out by a system of planes on $Q$. Supposing $\sD$ is rational over $K$,
the system of planes must be as well (but note that $|\sD|$ itself does not need
to contain rational divisors, nor does the system need to contain any planes
rational over $K$). Using the parametrisation above, there is an
$(x:1)\in\PP^1(K)$ such that
$$Q^-=r^2+2x rs+ x^2 s^2.$$
It follows that $(0:0:0:x:-1)$ is the singular point of $Q$ and that $Q$
is a cone over the quadric in $\PP^3$ given by
$$Q^+(u,v,w)-r^2=0.$$
According to Lemma~\ref{lemma:quadsplit}, this quadric has two rational systems
of lines (and therefore, a cone over it has two rational systems of planes)
if
$$-\det Q^+=-\det(Q_1^++2xQ_2^++x^2Q_3^+)$$
is a square in $K$. It follows that
$F$ is isomorphic to
$$F: y^2=-\det(Q_1^++2xQ_2^++x^2Q_3^+).$$
The map $F\times F\to\Pic^0(D)$ described above gives rise to an isomorphism
$$\Prym(D/C)\simeq \Jac(F).$$

In Cases 2 and 3 above, the Jacobian of $D$ contains elliptic curves $E_1$,
$E_2$. In these cases, $D$ is in fact a double cover of genus $1$ curves
$C_1$ and $C_2$ with $\Jac(C_i)\cong E_i$. The
cover can be constructed explicitly in the following way.
The $Q\in L_i$ have a fixed singularity on the line $\{u=v=w=0\}$.
Projecting from that singularity yields an intersection of two quadrics in
$\PP^3$, a model of $C_i$. Let $\sigma_i\in\Aut_\Kbar(D)$ denote the involution
of $D$ over $E_i$. We have $\iota=\sigma_1\sigma_2$. 
In fact, the projection $(u:v:w:r:s)\to(u:v:w)$ corresponds
to $D\to D/\langle\sigma_1,\sigma_2\rangle$. This shows that the
canonical model of $C$ is in fact of genus $0$ and hence that $C$ is
hyperelliptic. This places us in the situation of Section~\ref{sec:prymhyp}.

\begin{table}
\begin{center}
\begin{tabular}{c|c|c|c|c}
Case&$C$ & $D$ & $\Gamma^-$ & $\Prym(D/C)$\\[0.2em]
\hline&&&\\[-0.8em]
1&Hyperelliptic & Hyperelliptic & --- & $\Jac(F)$\\[0.2em]
2&Hyperelliptic & Non-hyperelliptic & split singular & $E_1\times E_2$\\[0.2em]
3&Hyperelliptic & Non-hyperelliptic & non-split singular &
     $\Re_{K(\sqrt{d})/K}(E)$\\[0.2em]
4&Non-hyperelliptic & Non-hyperelliptic & non-singular & $\Jac(F)$
\end{tabular}

\medskip
\caption{Structure of the Prym variety}\label{tab:prymstruct}

\begin{tabular}{c|l}
Case&Models\\
\hline
1&
$\begin{array}{rl}
C:& y^2=Q(x)R(x)\mbox{ where }\deg(Q)=2, \deg(F)=6\\
D:& y_1^2=Q(x) \mbox{ and } y_2^2=R(x)\\
F:& y_2^2=R(x).
\end{array}$\\
\hline
2&
$\begin{array}{rl}
C:& y^2=R_1(x)R_2(x)\mbox{ where }\deg(R_1)=\deg(R_2)=4\\
D:& y_1^2=R_1(x) \mbox{ and } y_2^2=R_2(x)\\
E_1:& \Jac(y_1=R_1(x))\\
E_2:& \Jac(y_2=R_2(x)).
\end{array}$\\
\hline
3&
$\begin{array}{rl}
C:& y^2=N_{K(\sqrt{d})[x]/K[x]} R(x)\mbox{ where }\deg(R)=4\\
D:& (x,y_0,y_1)\mbox{ satisfying }(y_0+y_1\sqrt{d})^2=R(x)\\
E:& \Jac(y=R(x))
\end{array}$\\
\hline
4&
$\begin{array}{rl}
C:&Q_1(u,v,w) Q_3(u,v,w) = Q_2(u,v,w)^2 \mbox{ with $Q_i$ quadratic forms}\\
D:&\left\{\begin{array}{rcl}
 Q_1(u,v,w)&=&r^2\\
 Q_2(u,v,w)&=&rs\\
 Q_3(u,v,w)&=&s^2
 \end{array}\right.\\
F:&y^2=-\det(Q_1+2xQ_2+x^2Q_3)
\end{array}$
\end{tabular}
\end{center}

\medskip
\caption{Models of involved curves}
\label{tab:models}
\end{table}

\begin{theorem}\label{thrm:prymclass}
Let $K$ be a field of characteristic $0$ and let $C$ be a curve of genus $3$
over $K$ with an unramified double cover $D/C$. Then $\Prym(D/C)$ can be described as
given in Table~\ref{tab:prymstruct} depending on the nature of $C$ and $D$.
Models of the curves involved can be described as in Table~\ref{tab:models}. For
hyperelliptic curves, it is assumed they are hyperelliptic over a $\PP^1$ over
$K$.
\end{theorem}

\section{Mapping $D$ into $\Prym(D/C)$}
\label{sec:embD}

As was shown in Section~\ref{sec:prym}, if $C$ is hyperelliptic then $D$ is a
cover of the curves that span $\Prym(D/C)$. Hence, it is obvious how to map $D$
into $\Prym(D/C)$. In this section we show how $D$ can be mapped into
$\Prym(D/C)$ if $C$ is non-hyperelliptic.

First, if we have a rational point $P_0\in D(K)$, we can embed $D$ in
$\Jac(D)$ via the \emph{Abel-Jacobi map}
$$\begin{array}{ccc}
D&\to&\Jac(D)\\
P&\mapsto&[P-P_0].
\end{array}$$
When we combine this map with the projection map
$(\id_*-\iota_*):\Jac(D)\to\Prym(D/C)$,
we obtain the \emph{Abel-Prym map}
$$\begin{array}{ccc}
D&\to&\Prym(D/C)\\
P&\mapsto&[P-\iota(P)]-[P_0-\iota(P_0)].
\end{array}$$

In general, we do not have a rational base point $P_0$ at our disposal. We give
an alternative map, based on the description of $\Prym(D/C)$ as $\Jac(F)$ for
some component $F$ of $W^1_4(D)$.

Let $P_0$ be a point on $D$ and let $L$ be the tangent line of $C$ at $\pi(P)$.
Then $L\cdot C$, being a linear section of a canonical model, determines an
effective canonical divisor on $C$. Consequently,
$\pi^*(L\cdot C)=2P_0+2\iota(P_0)+P_1+\iota(P_1)+P_2+\iota(P_2)$ is an
effective canonical divisor on $D$.

We use the notation from
Section~\ref{sec:genus5}. Furthermore, we write $\tng_D(P)$ for the tangent
line of $D$ at $P$ and for a quadric $Q\subset\PP^n$ and $P,T\in\PP^n$ we write
$P^t Q T$ for the matrix product, where $Q$ is identified with its representing
$(n+1)\times(n+1)$ symmetric matrix and $P,T$ are interpreted as
$(n+1)$-dimensional column vectors of projective coordinates.

\begin{lemma}\label{lemma:phimap} With the notation as above, let $P$ be a point
in $D(\Kbar)$.

\begin{itemize}
\item[(i)] There are two divisors $\fD_1,\fD_2$ (with $\fD_1=\fD_2$ in
degenerate cases) such that $\fD_i\geq 2P$ and $[\fD_i]\in F(\Kbar)\subset
W^1_4(\Kbar)$.
\item[(ii)] If $\fD_1=2P+P_1+P_2$ then $\fD_2=2P+\iota(P_1)+\iota(P_2)$.
\item[(iii)] Let $T\in\tng_D(P)$ with $T\neq P$. Then
$x_*(\fD_i)$ satisfies
$$(T^t\cdot Q_1\cdot T)+2x(T^t\cdot Q_2\cdot T)+x^2(T^t\cdot Q_3\cdot T)=0.$$
\item[(iv)] The map
$$\begin{array}{ccccc}
\phi:&D(\Kbar)&\to&\Div^2(F/\Kbar)\\
&P&\mapsto&[\fD_1]+[\fD_2]
\end{array}$$
is defined over $K$.
\end{itemize}
\end{lemma}

\begin{proof}
(iii): By Lemma~\ref{lemma:quadplanes} we have that if $\fD\in F$ with
$\fD=2P+P_1+P_2$, then $\tng_D(2P)\subset V_\fD\subset Q_\fD$. This is exactly
the case if $T\in Q_\fD$ for some (and hence for all) $T\in \tng_D(P)$ with
$T\neq P$. Using that $\Gamma^-\simeq\PP^1$ via $(1:2x:x^2)\mapsto (x:1)$, we
get the desired equation.

\medskip\noindent
(i):  From (iii), we know that there are two quadrics $Q\in \Gamma^-$
containing $\tng_D(P)$. We show that such a conic is $Q_\fD$ for some $\fD \in
F$ such that $\fD\geq 2P$. Consider the plane $V$ spanned by $\tng_D(P)$ and a
singular point $(0:0:0:x:-1)$ of $Q$. Then $V\subset Q$ and by
Lemma~\ref{lemma:quadplanes}, $\fD=V\cdot D$ represents a divisor class in
$W^1_4(D)$ and $Q=Q_\fD$.

If there is another divisor $\fD'\geq 2P$ with $Q_{\fD'}=Q$, then
$\fD-\fD'$ or $\fD-\iota_*(\fD')$ is principal. Taking the image under $\pi_*$
would give a divisor of a degree $2$ function on $C$, which contradicts that $C$
is not hyperelliptic.

\medskip\noindent
(ii): Note that $\pi_*(\fD_2)=\pi_*(\fD_1)=C\cdot \tng_C(\pi(P))$. If
$\fD_2=2P+P_1+\iota(P_2)$ then $2P+P_1$ must lie on the line
$V_{\fD_1}\cap V_{\fD_2}$. Since $D$ is canonical, it follows that
$l(2P+P_1)=2$. This contradicts Lemma~\ref{lemma:nontrig}. It follows that
$\fD_2$ must be as stated.

\medskip\noindent
(iv): Verify that $\act{\sigma}{(\phi(P))}=\phi(\act{\sigma}{P})$ for
$\sigma\in\Gal(\Kbar/K)$ via direct computation.
\end{proof}

\begin{prop}\label{prop:Demb} Let $C$ be a non-hyperelliptic genus $3$ curve
over a field $K$ of characteristic $0$ and let $D/C$ be an unramified double
cover with $\iota:D\to D$ the associated involution. Let $F$ be the genus $2$
curve given by Theorem~\ref{thrm:prymclass} such that $\Jac(F)=\Prym(D/C)$
and let $\phi:D\to\Div^2(F)$ be the map defined in Lemma~\ref{lemma:phimap}.
Then we have

$$\begin{array}{cccc}
2(\iota_*-\id_*):&D(\Kbar)&\to&\Prym(D/C)(\Kbar)\\
                 &P&\mapsto&\ [\phi(P)-\kappa_F].
\end{array}$$
\end{prop}
\begin{proof}
Using (ii) of Lemma~\ref{lemma:phimap}, we have
$\phi(P)=4P+P_1+P_2+\iota(P_1)+\iota(P_2)$. Using that
$[\kappa_D]=[2P+2\iota(P)+P_1+P_2+\iota(P_1)+\iota(P_2)]$, we have
$[\phi(P)-\kappa_D]=2P-2\iota(P)$ as an element of $\Pic(D)$.
Note that $[\kappa_D]=[\fD+\iota(\fD)]$ for
any $[\fD]\in F\subset W^1_4(D)$, so identifying
$\Pic(F)\subset\Pic(D)$, we get $[\kappa_F]=[\kappa_D]$. This proves the
proposition.
\end{proof}

\begin{remark}
It is worth noting that the map $(\iota_*-\id_*):D\to\ker(\pi_*)$ does not map
$D$ to $\Prym(D/C)$. To see this, note that $\Prym(D/C)=\Jac(F)$ for some genus
$2$ curve $F$. Any degree $0$ divisor class on $F$ can be represented as the
difference of two points on $F$, so we have that $\Prym(D/C)=F-F$. If  $P\in
D(\Kbar)$ has $[\iota(P)-P]\in \Prym(D/C)$, then we can find $[\fD_1],[\fD_2]\in
F(\Kbar)$ such that $[\iota(P)-P]=[\fD_1-\fD_2]$. Since $F\subset W^1_4(D)$,
we can choose effective representatives $\fD_1,\fD_2$
with a given point in the support. Hence, we can assume that
$\fD_1=\iota(P)+P_2+P_3+P_4$ and $\fD_2=P+P_5+P_6+P_7$. It follows that
$0=[P_2+P_3+P_4-P_5-P_6-P_7]$, so $\iota$ has a fixed point or $D$ is
hyperelliptic or trigonal. Our assumptions and Lemma~\ref{lemma:nontrig} tell us
that none of these are the case.
\end{remark}

Lemma~\ref{lemma:phimap} together with Proposition~\ref{prop:Demb} provide an
explicit way of mapping $D$ into an Abelian surface $\Prym(D/C)\simeq\Jac(F)$.
By slight abuse of notation, we write $\phi:D\to\Jac(F)$.

For explicit computations, Abelian surfaces have proven to be rather unwieldy.
In many cases, enough of the variety structure remains in the associated
Kummer-surface $\sK=\Jac(F)/\langle -1 \rangle$. The surface $\sK$ is naturally
expressed as a quartic surface in $\PP^3$. The map
$$\begin{array}{ccccc}
k:&\Jac(F)&\to&\sK\\
&[(x_1,y_1)+(x_2,y_2)-\kappa_F]&\mapsto&(1:x_1+x_2:x_1x_2:\ldots)&
                                                  =(k_1:k_2:k_3:k_4)
\end{array}$$
expresses $\Jac(F)$ as a double cover of $\sK$, ramified at $\Jac(F)[2]$, which
maps to the singular locus of $\sK$. The equation of $\sK$ is of the form
$$\sK: (k_2^2-4k_1k_3) k_4^2+K_3(k_1,k_2,k_3) k_4+ K_4(k_1,k_2,k_3)=0,$$
where $K_3,K_4$ are homogeneous forms of degrees $3,4$ respectively (see
\cite{prolegom} for explicit formulae). Hence,
$\sK$ is itself a double cover of the projective plane with coordinates
$(k_1:k_2:k_3)$ outside the point $(0:0:0:1)$.

Since $\iota_*\circ 2(\iota_*-\id_*)=-1\circ 2(\iota_*-\id_*)$, we see that
$D\to k\phi(D)$ factors through $D/\langle\iota\rangle=C$.
Furthermore, if $\fD\in \Div^2(F)$ is effective and
$(k_1:k_2:k_3:k_4)=k([\fD-\kappa_F])$, then $x_*(\fD)$ satisfies
$$k_3-k_2 x+k_1 x^2=0.$$
This gives us a procedure to compute many pointwise images for $k\phi$:
\begin{enumerate}
\item Choose an extension $L$ and a point $P\in D(L)$ (since $D$ is given as a degree
$8$ curve, there is an abundance of suitable degree $8$ extensions)
\item Following Lemma~\ref{lemma:phimap}, choose $T\in\tng_D(P)$ and set
$$(k_1,k_2,k_3)=
   (T^t\cdot Q_1 \cdot T,-2T^t\cdot Q_2 \cdot T,T^t\cdot Q_3 \cdot T)$$
\item If $k_3-k_2x+k_1 x^2$ is irreducible of degree $2$ over $L$, then there is
a unique point $(k_1:k_2:k_3:k_4)\in\sK(L)$ that has an $L$-rational preimage in
$\Jac(F)$. This is the desired image.
\end{enumerate}
The irreducibility in the last step corresponds to $P_1,P_2$ from
Lemma~\ref{lemma:phimap}(ii) being quadratic conjugate over $L$. In that case,
the divisor $P_1+\iota(P_2)$ is not $L$-rational and hence rationality tells
which divisor to pick. If $P_1,P_2$ are themselves $L$-rational, then the above
procedure does not compute sufficient information to distinguish between
$2P+P_1+P_2$ and $P+\iota(P)+P_1+\iota(P_2)$.

The procedure above yields a way to compute the equations of $k\phi(D)$. First,
one gathers many pointwise images for points over extensions $L$ and then one
interpolates for low degree rational forms vanishing on those points. As we will
see, $k\phi(D)$ is the intersection of $\sK$ with another degree $4$ surface.

\begin{lemma} With the notation above, the image of $D$ under
$k\circ\phi$ is of degree at most $16$.
\end{lemma}

\begin{proof}
We compute the degree of the image by computing the degree of the intersection
with $k_1=0$. By change of basis we can assume that $Q_3$ is of rank $4$.
A point $P\in D$  has $k_1=0$ if $Q_3$ contains the plane $V_\fD$
for some $\fD\geq 2P$. The two plane rulings on $Q_3$ give rise to two degree $4$
covers $D\to\PP^1$, where the fibres are the divisors cut out by the $V_\fD$ in
the ruling. From Riemann-Hurwitz it follows that for each ruling there are $16$
ramified fibres, i.e., $\fD$ of the form $2P+P_1+P_2$. Hence, we see that there
are $32$ points on $D$ (counted with appropriate multiplicity) that land on
$k_1=0$. Note that the image of $D$ under $k\phi$ factors through
$D/\langle\iota\rangle$ , so the degree of $k\phi(D)$ is at most $16$.
\end{proof}

The procedure above has been implemented as a routine for the computer algebra
system MAGMA \cite{magma}. See \cite{bruin:prymelectr}.

\section{The fibre of the Prym map in genus $3$}\label{sec:fibprym}

Given a genus $2$ curve $F$, we have an Abelian variety $\Jac(F)$ and a quartic
surface $\Jac(F)/\langle -1 \rangle=\sK\subset\PP^3$, with a singular locus
consisting of the image of $\Jac(F)[2]$. There is an obvious way of realising
$\Jac(F)$ as a Prym variety of a non-hyperelliptic curve of genus $3$. Pick a plane $V\subset
\PP^3$ such that $C:=V\cap\sK$ is a non-singular quartic curve. It follows that
$C$ stays away from the singular points of $\sK$ and thus does not meet the
ramification locus of $k:\Jac(F)\to\sK$. Therefore $D=k^{-1}(C)$ is an
unramified cover of $C$. Either $D$ is connected and hence of genus $5$ or $D$
is the disjoint union of two copies of $C$. Note however that $C$ is of genus
$3$ and hence has to be special to fit in an Abelian surface.

In fact, as Verra~\cite{verra:prym} proves, over an algebraically closed field,
essentially any occurrence of $\Jac(F)$ as $\Prym(D/C)$ occurs for $C$
isomorphic to a linear section of $\sK$. The addition of $\Jac(F)[2]$ induces
automorphisms of $\sK$ which are induced by linear transformations of $\PP^3$.
He shows that the fibre of the Prym map $(D/C)\mapsto\Prym(D/C)$ over $\Jac(F)$
is a blow-up of $\widehat{\PP}^3/\Jac(F)[2]$, where $\widehat{\PP}^3$ is the
space of plane sections of $\sK$.

Verra also proves that genus $5$ curves $D\subset\Jac(F)$ of the form above
are, up to translation by a $2$-torsion point, Abel-Prym embeddings. We give a
short description of a procedure to recover 
$$C: Q_1(u,v,w)Q_3(u,v,w)=Q_2(u,v,w)^2$$
from a plane section of the Kummer surface $\sK=\Jac(F)/\langle -1\rangle$ for a
curve $F$ of genus $2$.

First we review some of the basic geometry of Jacobians of curves of genus $2$.
We follow the notation introduced in \cite{prolegom}. Let $F$ be a curve of
genus $2$. They define a projective model of $\Jac(F)$ in
$\PP^{15}$ with coordinates $(z_0:\cdots:z_{15})$ with (among others) the following properties:
\begin{itemize}
\item There is a symmetric theta-divisor $\Theta$ on $\Jac(F)$ such that
$$\langle k_1,\ldots,k_4\rangle=|2\Theta|$$
and
$$\langle z_0,\ldots,z_{15}\rangle=|4\Theta|$$
with
$$(k_1:k_2:k_3:k_4)=(z_{14}:z_{13}:z_{12}:z_5).$$
\item  The coordinates $(k_1:k_2:k_3:k_4)$ provide a model of the Kummer surface
$\sK=\Jac(F)/\langle -1\rangle$.
\item With respect to the action of $-1\in\Aut(\Jac(F))$ on $|4\Theta|$, we have
that 
$$\langle z_0,z_3,z_4,z_5,z_{10},\ldots,z_{15}\rangle=\langle
k_1^2,k_1k_2,\ldots,k_4^2\rangle$$
is the $+1$-eigenspace and
$$\langle g_0,\ldots,g_5\rangle:=\langle z_1,z_2,z_6,z_7,z_8,z_9\rangle$$
is the $-1$-eigenspace.
\end{itemize}

Let $V=\{k_4=v_1k_1+v_2k_2+v_3k_3\}\subset\PP^3$ be a plane such that $C=\sK\cap V$ is a non-singular plane
section and let $(u:v:w)=(k_1:k_2:k_3)$ be coordinates on $V$ (since $(0:0:0:1)$
is a singular point of $\sK$, assuming the suggested form of $V$ is not a
restriction).
 The curve $C$ is a nonsingular degree
$4$ plane curve. It follows that $C$ is of genus $3$ and that $(u:v:w)$ gives a
canonical model of $C$, i.e., that $\langle u|_C,v|_C,w|_C\rangle=|\kappa_C|$
for some canonical divisor $\kappa_C$ of $C$. Let $D$ be as above and let
$\kappa_D$ be the pull-back of $\kappa_C$. It is a straightforward computation
to check that the restriction of $\langle z_0,\ldots,z_{15} \rangle$ to $D$
gives a linear system contained in
$|\kappa_D|$ and that generically it gives the complete linear system.

Using the quadratic relations between the $z_i$ (see \cite{flynn:formgroup}, \cite{prolegom} and
\cite{flynn:electr} for the explicit formulae), we can express any $g_ig_j$ as
a degree $4$ form in $k_1,\ldots,k_4$. Hence, we obtain that
$$(a_0g_0+\cdots+a_5g_5)^2=G(a_0,\ldots,a_5;k_1,\ldots,k_4),$$
where $G$ is homogeneous of degrees $2,4$ in the $a_i$ and the $k_j$
respectively.

Insisting that
$$G(a_0,\ldots,a_5)=u^2 Q(u,v,w)$$
for some quadratic form $Q$ gives $9$ quadratic equations in $a_0,\ldots,a_5$.
However, a solution to these equations corresponds exactly to a form
$Q(u,v,w)$ on $C$ that becomes a square when pulled back to $D$. We know that
this happens for exactly two forms $Q_1^+(u,v,w)=r^2$ and $Q_3^+(u,v,w)=s^2$, so
these equations determine a degree $2$ locus in $(a_0:\cdots:a_5)$.

Solving these equations allows us to determine $Q_1^+(u,v,w)$ and $Q_3^+(u,v,w)$
up to a scalar. The quadratic form $Q_2^+(u,v,w)$ is then easily determined up
to a scalar because this corresponds to the conic through the intersection of
$Q_1^+(u,v,w)Q_3^+(u,v,w)=0$ with $C$. The appropriate scalars $\lambda,\mu$ are
then easily determined from the fact that
$$\lambda Q_1^+(u,v,w)\mu Q_3^+(u,v,w)-(\mu Q_2^+(u,v,w))^2$$
should equal the equation for $C$.

Note that we may find $Q_i^+$ that are quadratic
conjugate over the base field we are working with, because of the arbitrary
coordinate choice we made when insisting that $G=u^2 Q(u,v,w)$. The analysis
from
Section~\ref{sec:genus5} guarantees us that by change of basis (corresponding to
$\Aut(\Gamma^-)$ or, equivalently, fractional linear transformations of $x$) we
can in fact obtain rational $Q_1^+, Q_3^+$. This yields the following
amusing result, which is equivalent to saying that all Jacobians of genus $2$
curves over $\QQ$ occur as Prym varieties of non-hyperelliptic curves of genus
$3$ over $\QQ$.

\medskip
\noindent\textit{Proof of Proposition~\ref{prop:fdet}}:
Take the curve of genus $2$
$$F:y^2=f(x)$$
and take a sufficiently general plane section of the associated Kummer-surface.
Using the construction above, we obtain a cover $D\to C$ such that
$\Prym(D/C)=\Jac(F)$. Section~\ref{sec:genus5} tells us that $F$ must be of the
described form and the outline above explains how one can find the
representation explicitly.
\hfill$\square$
\medskip

In fact, the model of $C$ as a plane section of a Kummer-surface $\sK$
completely encodes the 28 bitangents of $C$ as well. The $16$ tropes of $\sK$
obviously cut out bitangents on $C$. The remaining $12$ bitangents come in
pairs, making up the $6$ singular conics in the family
$Q_1^+ +2x Q_2^+ +x^2 Q_3^2$.

This gives us a way to search for genus $3$ curves with all bitangents
rational. First, start with a Kummer-surface $\sK$ with $16$ rational tropes (i.e.,
the Kummer-surface of the Jacobian of a genus $2$ curve with $6$ rational
{Weierstra\ss} points). Then, select a plane $V$ such that $C:=V\cap \sK$ is of
the form
$$Q_1^+(u,v,w)Q_3^+(u,v,w)=Q_2^+(u,v,w)^2,$$
where the singular conics in $Q_1^+ +2x Q_2^+ + x^2 Q_3^2$ are split. As it turns
out, these are all split or nonsplit simultaneously.

\medskip
\noindent\textbf{Example:} \textit{A curve of genus $3$ with all $28$ bitangents
rational.}

\medskip
\noindent Take
$$F:y^2=x(x-2)(x-1)(x+1)(x+3).$$
The corresponding Kummer-surface is
$$\begin{array}{l}
\sK:36k_1^4 + 84k_1^3k_3 - 24k_1^2k_2k_3 - 12k_1^2k_2k_4 + 65k_1^2k_3^2 +
        4k_1^2k_3k_4-24k_1k_2^2k_3 +\\
\;\;4k_1k_2k_3^2 +14k_1k_2k_3k_4 + 14k_1k_3^3 - 4k_1k_3^2k_4 - 4k_1k_3k_4^2 +
    k_2^2k_4^2 - 2k_2k_3^2k_4 + k_3^4=0.
 \end{array}$$
We take the plane section
$$V: k_1+k_2+t k_3+k_4=0.$$
Projecting onto $(u:v:w)=(k_1:k_2:k_3)$, we obtain
$$C:Q_1^+(u,v,w) Q_3^+(u,v,w) = Q_2^+(u,v,w)^2$$
where
$$\begin{array}{rcl}
Q_1^+&=&(36t-82)u^2+(6t-80)uv+(-11t+14)uw+(6t+2)v^2+\\
    &&(t+14)vw+tw^2,\\
Q_2^+&=&(-3t+40)u^2+(-\frac{1}{2}t-9)uv+(-6t^2+\frac{11}{2}t+51)uw+\\
    &&(-\frac{1}{2}t-7)v^2+(-\frac{1}{2}t^2-2t+2)vw+(-t^2+\frac{1}{2}t+7)w^2,\\
Q_3^+&=&(6t+2)u^2+(t+14)uv+(t^2+2t-4)uw+tv^2+(2t^2-t-14)vw+\\
    &&(t^3-t^2-8t+2)w^2.
\end{array}$$

%$$\begin{array}{rcl}
%Q_1^+&=&64906(u + \frac{5}{83}v - \frac{13}{83}w)(u + \frac{11}{17}v -
%\frac{1}{17}w)\\
%Q_2^+&=&-22954u^2 + 4462uv - 19634uw + 3404v^2 - 1916vw -
%2728w^2\\
%Q_3^+&=&2530(u - \frac{13}{5}v + \frac{223}{115}w)
%         (u - \frac{1}{11}v - \frac{227}{253}w).
%\end{array}$$
The $16$ bitangents coming from the tropes are given by the polynomials
$$\begin{array}{l}
    u,\; w,\; 4u - 2v + w,\;
    5u + 3v + (t + 1)w,\;
    4u + v + (t + 2)w,\;
    u + 7v + (t - 1)w,\;\\
    u - v + w,\;
    7u + v + (-t - 3)w,\;
    9u + 3v + w,\;
    7u + v + (t + 1)w,\;\\
    5u - v + (-t + 1)w,\;
    u - v + (-t + 3)w,\;
    2u - 4v + (-t + 2)w,\;\\
    10u - 2v + (t - 4)w,\;
    u + v + w,\;
    2u + 2v + tw.
\end{array}$$
%$$\begin{array}{l}u - \frac{1}{5}v + \frac{49}{115}w,\;
%    u - v + \frac{95}{23}w,\;
%    w,\;
%    u - \frac{1}{2}v + \frac{1}{4}w,\;
%    u + \frac{3}{5}v - \frac{3}{115}w,\;
%%    u + \frac{1}{4}v + \frac{5}{23}w,\;
%    u,\;\\
%    u + v + w,\;
%    u - \frac{1}{5}v - \frac{59}{115}w,\;
%    u - 2v + \frac{36}{23}w,\;
%    u - v + w,\;
%    u + \frac{1}{3}v + \frac{1}{9}w,\;
%    u + \frac{1}{7}v - \frac{43}{161}w,\;\\
%    u + 7v - \frac{49}{23}w,\;
%    u + v - \frac{13}{23}w,\;
%    u + \frac{1}{7}v - \frac{3}{161}w.
% \end{array}$$
The remaining $12$ bitangents come from the $6$ singular quadrics. They are
split if $196+20t-23t^2$ is a square. Therefore, we substitute
$$t:=\frac{4s^2-10s-6}{2s^2+s+3}$$
and obtain the bitangents given by the polynomials
$$\begin{array}{l}
(s+21)u+(7s+3)v+(s-3)w,\;
(10s+11)u+(-2s+5)v+(-2s-1)w,\\
(8s^3+2s^2+11s-3)u+(2s^3+5s^2+5s+6)v+(8s^3+17s^2-s-6)w,\\
(10s^3-s^2+12s-9)u+(-2s^3-7s^2-6s-9)v+(-2s^3+5s^2+30s-9)w,\\
(4s^3+4s^2+7s+3)u+(2s^2+s+3)v+(2s^2+3s-3)w,\\
(2s^3-5s^2-9)u+(-2s^3-3s^2-4s-3)v+(2s^3+7s^2+4s-1)w,\\
(2s^3+7s^2+6s+9)u+(2s^3+s^2+3s)v+(2s^3+s^2-3s)w,\\
(4s^3+10s^2+10s+12)u+(-2s^3+s^2-2s+3)v+(s^3+4s^2-5s)w,\\
(4s^3+16s^2+13s+21)u+(-4s^3-5s+3)v+(4s^3+16s^2-3s-3)w,\\
(10s^3+23s^2+24s+27)u+(6s^3+s^2+8s-3)v+(6s^3-5s^2-20s+7)w,\\
(14s^3+13s^2+24s+9)u+(2s^3-5s^2-9)v+(-10s^3-35s^2+9)w,\\
(4s^3-8s^2+s-15)u+(4s^3+4s^2+7s+3)v+(4s^3-8s^2-23s+9)w.
\end{array}$$

%$$\begin{array}{l}
%u - \frac{1}{11}v - \frac{227}{253}w,\;
%    u - \frac{13}{5}v + \frac{223}{115}w,\;
%    u + \frac{11}{17}v - \frac{1}{17}w,\;
%    u + \frac{5}{83}v - \frac{13}{83}w,\;
%    u + \frac{11}{17}v + \frac{517}{391}w,\;\\
%    u - \frac{7}{8}v + \frac{77}{92}w,\;
%    u + 2v - \frac{58}{23}w,\;
%    u + \frac{3}{13}v + \frac{51}{299}w,\;
%    u + \frac{5}{14}v + \frac{95}{644}w,\;
%    u - \frac{1}{11}v + \frac{25}{253}w,\;\\
%    u + \frac{3}{13}v - \frac{21}{299}w,\;
%    u - \frac{7}{31}v + \frac{371}{713}w.
%\end{array}$$

\section{Applications to finding rational points on curves of genus $3$}

In this section, we will apply the concepts of \emph{covering collections} (see
\cite{chevweil}, \cite{weth:phdthesis}, \cite{bruin:tract}, \cite{brufly:tow2cov}) and
\emph{Chabauty methods} (see \cite{cole:effchab}, \cite{flynn:flexchab}) to a
curve $C$ of genus $3$ with an unramified double cover $D$. We end up
determining the rational points on a curve of genus $5$ inside the
Jacobian of a curve $F$ of genus $2$. The hardest piece of information we will
need is the Mordell-Weil group of $\Jac(F)$. Computationally, this is much more
attractive than applying Chabauty methods directly to an embedding of $C$ in
its own Jacobian. In the latter case, we would have to analyse the Mordell-Weil
group of $\Jac(C)$.

Additionally, the techniques we present here do not depend on the existence of
an embedding of $C$ in $\Jac(C)$. As a result, we will see that we can even use
the construction to exhibit part of a local-global obstruction for $C$ and $D$
having rational points.

Let $K$ be a number field and let $C$ be a non-hyperelliptic curve of genus $3$
with an unramified double cover $D$ over $K$. As we have seen in
Section~\ref{sec:genus5}, it follows that there exists a smooth plane model of
$C$ of the form
$$Q_1(u,v,w)Q_3(u,v,w)=Q_2(u,v,w)^2,$$
where $Q_1,Q_2,Q_3\in K[u,v,w]$ are quadratic forms. Without loss of generality,
we can assume that $Q_1,Q_2,Q_3$ have integral coefficients. Furthermore, we
have a collection of twists of $D$, each covering $C$, of the form:
$$D_\delta:\left\{
\begin{array}{rcl}
Q_1(u,v,w)=\delta r^2\\
Q_2(u,v,w)=\delta rs\\
Q_3(u,v,w)=\delta s^2\\
\end{array}
\right.$$
We write $\sO=\sO_K$ for the ring of integer of $K$ and
we consider the projective $\sO_K$-scheme $X$ corresponding to the ideal
$I=(Q_1,Q_3,Q_1Q_3-Q_2^2)\sO_K[u,v,w]$. Since $C$ is non-singular as a curve
over $K$, we have that $X\times_{\sO_K}\Spec(K)$ is empty. Let $S$ be a finite set of primes
such that $X\times\Spec(\sO_S)$ is empty. Such a set $S$ is easily computed. For
instance, compute
$$\Res_u(\Res_v(Q_1,Q_3),\Res_v(Q_1,Q_2))=\lambda w^{16}.$$
One can take $S$ to be the set of prime divisors of $\lambda$. One may obtain a
smaller set by intersecting such sets $S$ obtained from all different
combinations in which such resultants could be taken.
We recall the definition
$$K(S,2):=\{\delta\in K^*: v_\fp(\delta)\equiv 0\mod 2\mbox{ for all primes
of $K$ satisfying } \fp\notin S\}/ K^{*2}.$$
This is a finite subgroup of $K^*/K^{*2}$ and we will identify its elements with
a finite set of representatives in $K^*$.

We obtain the standard lemma:
\begin{lemma} Let $C$ and $Q_1,Q_2,Q_3$ and $S$ be as above. If $(u_0:v_0:w_0)\in C(K)$,
then there exists $\delta\in K(S,2)$ and $r_0,s_0\in K$ such that
$$(u_0:v_0:w_0:r_0:s_0)\in D_\delta(K).$$
\end{lemma}

It follows that any rational point on $C$ has a rational pre-image on $D_\delta$
for some $\delta \in K(S,2)$. Thus, in order to determine the rational points of
$C$, it suffices to determine the rational points of $D_\delta$ for all $\delta\in
K(S,2)$. From Section~\ref{sec:W14} we know that for
$$F_\delta: y^2=-\delta \det(Q_1+2x Q_2 + x^2 Q_3),$$ 
we have $\Prym(D_\delta/C)\simeq\Jac(F_\delta)$ and Proposition~\ref{prop:Demb}
gives an explicitly computable map $\phi:D_\delta\to\Jac(F_\delta)$.
We can then proceed to determine $\phi(D_\delta(\QQ))\cap\Jac(F_\delta)(\QQ)$
or rather, as it turns out, $k(\Jac(F_\delta)(\QQ))\cap k\phi(D_\delta)(\QQ)$.

\medskip
\noindent\textbf{Example:} \textit{Chabauty using Prym varieties}.

\medskip
\noindent\emph{Proof of Proposition~\ref{prop:ex1}}: See
\cite{bruin:prymelectr} for a transcript of the computer calculations.
Applying the method described above we verify that we can take $S=\{1,2,5\}$ and local
considerations show that $D_\delta(\QQ)=\emptyset$ for $\delta\neq -1$. We find
$$F: y^2=x^5 + 8x^4 -  7x^3 -\frac{7}{2}x^2 + 5x - 1$$
and
$$\Jac(F)(\QQ)=\langle \sD \rangle=
\langle[(2\sqrt{2}-2,17\sqrt{2}-25)+(-2\sqrt{2}-2,-17\sqrt{2}-25)-2\infty]\rangle.$$
The equation of the associated Kummer-surface is
$$\begin{small}\begin{array}{ll}
\sK:&
11k_1^4 - 28k_1^3k_2 + 70k_1^3k_3 + 4k_1^3k_4 + 32k_1^2k_2^2 - 164k_1^2k_2k_3 -
10k_1^2k_2k_4 + 171k_1^2k_3^2 + \\
&14k_1^2k_3k_4 + 4k_1k_2^3 - 20k_1k_2^2k_3 + 14k_1k_2k_3^2 + 14k_1k_2k_3k_4 +
 14k_1k_3^3 - 32k_1k_3^2k_4 - \\
&4k_1k_3k_4^2 + k_2^2k_4^2 - 2k_2k_3^2k_4 + k_3^4=0
 \end{array}\end{small}$$
and, using the interpolation procedure described in Section~\ref{sec:embD}, we
find that the embedding of $C$ in $\sK$ as $k\phi(D)$ is given by the equation
$$\begin{tiny}\begin{array}{ll}
\psi:&
429136k_1^4+1330784k_1^3k_3+567232k_1^3k_4-159200k_1^2k_2^2-2866016k_1^2k_2k_3+
33440k_1^2k_2k_4+4248768k_1^2k_3^2+\\
&27552k_1^2k_3k_4+881664k_1^2k_4^2+288072k_1k_2^3-777432k_1k_2^2k_3-
256928k_1k_2^2k_4+244832k_1k_2k_3^2+\\
&907424k_1k_2k_3k_4-745472k_1k_2k_4^2+593152k_1k_3^3-991488k_1k_3^2k_4+
357440k_1k_3k_4^2+573440k_1k_4^3+34895k_2^4-\\
&69720k_2^3k_3+1120k_2^3k_4+151704k_2^2k_3^2-364448k_2^2k_3k_4+
226032k_2^2k_4^2-251552k_2k_3^3+569376k_2k_3^2k_4+\\
&10752k_2k_3k_4^2-315392k_2k_4^3+156704k_3^4-167552k_3^3k_4
-283136k_3^2k_4^2+200704k_3k_4^3+114688k_4^4=0.
\end{array}\end{tiny}$$
It is a straightforward computation to check that the intersection of $\sK$ with
$\psi(k)=0$ is non-singular, which verifies that $k\phi(D)$ is indeed an
embedding of $C$ in $\sK$ and that $\phi(D)$ is an embedding of $D$ in
$\Jac(F)$.

Using that $\phi:D\to\Jac(F)$ is defined over $\QQ$ and that
$\Jac(F)(\QQ)=\langle \sD\rangle$, we find that a rational point
$P\in\phi(D)(\QQ)$ must be of the form $P=n\sD$ for some $n\in \ZZ$.
Furthermore, considering the $F$ and $\sK$ over $\FF_{13}$, we find that any
such point must have $n\equiv \pm 1 \mod 10$.

Using the formal group law of $\Jac(F)$, we obtain a power series
$$\psi(N)=\psi(k((1+10N)\sD))=\psi(K(\sD+\Exp(N\Log(10\sD)))),$$
say,
$$\psi(N)=\psi_0+\psi_1 N+\psi_2 N^2+\cdots\in\ZZ_{13}[\![ N]\!],$$
with $\psi_i\equiv 0\mod 13^i$. 
such that, if $P=(1+10N)\sD$ is a point on $\phi(D)(K)$, then $\psi(N)=0$.

Note that the values of $\psi(k(\sD)),\psi(k(11\sD))$ determine
$\psi_0,\psi_1\mod 13^2$, so one does not need an explicit description of the
formal group law on $\Jac(F)$ to obtain an approximation to $\psi(N)$.

Since $\psi(k(\sD))=0$ and $\psi(k(11\sD))\not\equiv 0\mod 13^2$, it follows
that $\psi_1\not\equiv 0\mod 13^2$. From Stra{\ss}mann's lemma it follows that
$\psi(N)$ has at most one $0$ for $N\in\ZZ_{13}$ (i.e., $N=0$). This implies
that $\phi(D)$ has only one rational point which reduces to $\sD$ modulo $13$.
By symmetry, it follows that there is also only one rational point reducing to
$-\sD$ modulo $13$. On the other hand, the computation over $\FF_{13}$ shows
that all rational points of $\phi(D)$ reduce to $\pm\sD$ modulo $13$. Hence, it
follows that
$$\phi(D)(K)=\{\sD,-\sD\}$$
and that $C$ has only one rational point, being $(0:1:0)$.
\hfill$\square$
\medskip

\noindent\textbf{Example:} \textit{Computations in the Brauer-Manin obstruction}.

Since the embedding of $D$ in $\Prym(D/C)$ is independent of $D$ having any
rational points, we can also apply this construction to curves $D$ that have no
rational points, but do have rational points everywhere locally. Using
information obtained from the reduction of $\Jac(F)$ at various primes, we might
actually succeed in \emph{proving} that $D(\QQ)$ is empty. Under the assumption
that $\Jac(D)$ has a finite Tate-Shafarevich group, this corresponds to
computing part of the Brauer-Manin obstruction according to
\cite{schar:phdthesis}.

We consider the curve
$$C:(v^2 + vw - w^2)(uv + w^2) = (u^2 - v^2 - w^2)^2.$$
It is easily checked that $C$ has points everywhere locally. Furthermore, the
set of primes $S$ described above can be taken to be $\{2\}$ and only for
$\delta=1$ does $D_\delta$ have points everywhere locally.

\medskip
\noindent\emph{Proof of Proposition~\ref{prop:ex2}}: See
\cite{bruin:prymelectr} for a transcript of the computer calculations.
The fact that $D(\QQ_p)$ and $D(\RR)$ are all non-empty can be verified with a
straightforward computation. To prove that $D(\QQ)=\emptyset$ we embed $D$ in
$\Prym(D/C)=\Jac(F)$, where
$$F:y^2=x^6+2x^5+15x^4+40x^3-10x.$$
We find that
$$\Jac(F)(\QQ)=\langle \sD_1,\sD_2\rangle\simeq \ZZ\times \ZZ,$$
where
$$\begin{array}{rcl}
\sD_1&=&[\infty^+-\infty^-],\\
\sD_2&=&[\{x^2+\frac{2}{5}x+\frac{1}{41}=0,y=\frac{4683}{2050}x -
\frac{281}{410}\}\cdot F -\kappa_F].
\end{array}$$
Considering $\Jac(F)$ modulo $7$, we find
$$\phi(D)(\QQ)\subset\{\pm9\sD_1,\pm22\sD_1,\pm23\sD_1\}+\langle
55\sD_1,\sD_2-15\sD_1\rangle$$
and considering $\Jac(F)$ modulo $11$, we find
$$\phi(D)(\QQ)\subset\{\pm33\sD_1\}+\langle 93\sD_1,\sD_2+46\sD_1\rangle.$$
Considering $\Jac(F)$ modulo $11^2$, we find that the two residue classes modulo
$11$ lift to $11$ residue class modulo
$$\langle 11\cdot 93\sD_1,11\cdot(\sD_2+46\sD_1)\rangle.$$
We combine the information modulo $11^2$ and $7$ and express it as congruences
modulo
$$\langle 55\sD_1,\sD_2-15\sD_1,11\cdot 93\sD_1,11\cdot(\sD_2+46\sD_1)\rangle=
\langle 11\sD_1,\sD_2-4\sD_1\rangle.$$
We obtain
$$\begin{array}{lllllll}
\mbox{from
}7&:&\phi(D)(\QQ)&\subset&\{0,\pm\sD_1,\pm2\sD_1\}&+&\langle11\sD_1,\sD_2-4\sD_1\rangle\\
\mbox{from }11^2&:&\phi(D)(\QQ)&\subset&\{\pm 4\sD_1\}&+&\langle11\sD_1,\sD_2-4\sD_1\rangle
\end{array}.$$
It follows that $D(\QQ)=\emptyset$. and therefore that $C(\QQ)$ is empty as
well.
\hfill$\square$
\medskip

\section{Acknowledgements}

I would like to thank MSRI for their 2000 fall special semester on explicit
methods in number theory for bringing me into contact with many mathematicians.
I have greatly benefitted from the discussions I had with them. In particular I
would like to thank, Armand Brumer, J. Merriman, Ed Schaefer, Michael Stoll and
Joe Wetherell.

%\bibliography{art.bib}

\end{document}